\documentclass{article}
\usepackage{amsfonts,amsmath,amsthm,amssymb,graphics,epsfig}

\input{tcilatex}
\begin{document}

\title{Singular Points of Reducible Sextic Curves}
\author{David A. Weinberg \and Nicholas J. Willis}
\date{August 20, 2007}
\maketitle

\begin{abstract}
There are 106 individual types of singular points for reducible complex
sextic curves.
\end{abstract}

\section{Introduction}

\bigskip Extensive studies of simple singularities of complex sextic curves
have been made by T. Urabe [4][5] and J.-G. Yang [10]; see also [1]. \
Singular points of sextic curves of torus type have been classified by M.
Oka [2][3]. The authors have classified all individual types of singular
points for irreducible real and complex sextic curves in a previous paper
[8]. \ In this paper, we will determine the individual types of singular
points for reducible complex sextic curves. \ There are 106 types. \ The
proof is as elementary as possible and relies heavily on Puiseux expansions,
which are computed by using Maple when necessary.

All definitions pertaining to the classification, including the equivalence
relation itself, were explained in detail in previous papers [6], [7], [8].
\ Recall that we assign a diagram to each type of singular point and that
the assignment of a diagram is invariant with respect to linear changes of
coordinates.

Reducible sextic curves fall into families that must have an irreducible
factor of degree 1, 2, or 3. \ If there is an irreducible factor of degree
one, then every type of singular point in that family can be obtained by
careful scrutiny of the Newton polygon and a knowledge of all types of
singular points for quintic curves. \ Next consider the cases where there is
an irreducible factor of degree two or three (but none of degree one). \ If
this factor does not share a common tangent line with the remaining
factor(s), then the singular point types can be determined by mathematical
common sense. \ If this factor does share a common tangent line with any of
the remaining factor(s), then we studied the family carefully by means of
Maple computations.

All such families are considered in the next section. \ In the final section
of the paper, we give the complete list of 106 types of singular points for
reducible complex sextic curves.

\bigskip

Acknowledgment. \ \textit{The authors wish to thank Mark van Hoeij (Florida
State University) for some singularly valuable Maple code.}

\section{Results of Symbolic Computations}

\bigskip In this section we study the cases where the reducible sextic curve
has an irreducible factor of degree two or three (but none of degree one)
that shares a common tangent line with the remaining factor(s). \ Such
families were studied carefully by means of Maple computations. \ Since
these computations were quite lengthy, we refer the reader to the Maple
worksheets posted on the website of David Weinberg for the details [9]. \
Not every singular point of a sextic curve can be described by using the
traditional Arnol'd notation. \ In this paper, we will express the singular
point type by using the diagrams described in our previous papers [6], [7],
[8]. \ For each family, we indicate the diagrams giving the singular point
types.

\bigskip

1. $%
(y+ax^{2}+bxy+cy^{2})(y^{2}+dxy+ex^{3}+fx^{2}y+gxy^{2}+hy^{3}+jx^{4}+kx^{3}y+lx^{2}y^{2}+mxy^{3}+ny^{4})=0
$

conic with one branch through origin, quartic with two distinct branches
through origin, one tangent line in common
\begin{figure}[h]
  \caption{a=1; b=2, 3, 4, 5, 6, or 7.}
  \begin{center}
    \includegraphics{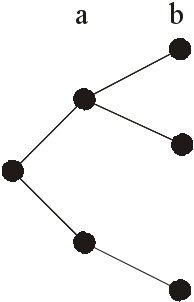}
  \end{center}
\end{figure}

2. \ $%
(y+ax^{2}+bxy+cy^{2})(y^{2}+ex^{3}+fx^{2}y+gxy^{2}+hy^{3}+jx^{4}+kx^{3}y+lx^{2}y^{2}+mxy^{3}+ny^{4})=0
$

conic with one branch thru origin, quartic with 2 branches through the
origin and one tangent line, the conic and the quartic share a common
tangent line.
\begin{figure}[h]
  \caption{a=3/2 or 2}
  \begin{center}
    \includegraphics{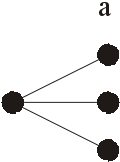}
  \end{center}
\end{figure}
\clearpage
\begin{figure}[h]
  \caption{a=2; b=5/2, 3, 4, 5, or 6}
  \begin{center}
    \includegraphics{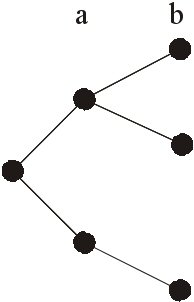}
  \end{center}
\end{figure}

 3. $%
(y+ax^{2}+bxy+cy^{2})(y^{3}+dx^{2}y+exy^{2}+fx^{4}+gx^{3}y+hx^{2}y^{2}+kxy^{3}+ly^{4})=0
$

conic with one branch thru origin, quartic with three branches thru origin,
one tangent line in common.

\begin{figure}[!h]
  \caption{a=1; b=2, 3, 4, 5, or 6}
  \begin{center}
{\includegraphics{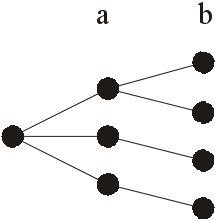}}
\end{center}
\end{figure}

\bigskip

\begin{figure}[!h]
  \caption{a=1; b=3/2; c=2}
  \begin{center}
  {\includegraphics{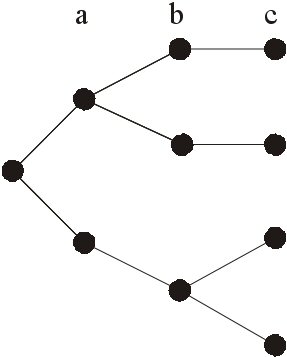}}
\end{center}
\end{figure}

 4. $\
(y+ax^{2}+bxy+cy^{2})(y+dx^{2}+exy+fy^{2}+gx^{3}+hx^{2}y+jxy^{2}+ky^{3}+lx^{4}+mx^{3}y+nx^{2}y^{2}+pxy^{3}+qy^{4})=0
$

conic with one branch thru the origin, quartic with one branch
thru origin, one tangent line in common.
\begin{figure}[!h]
  \caption{a=2, 3, 4, 5, 6, 7, or 8}
  \begin{center}
{\includegraphics{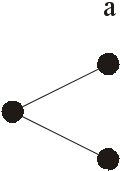}}
\end{center}
\end{figure}
\clearpage
5. \ $%
(y+ax^{2}+bxy+cy^{2})(y^{2}+dx^{2}y+ex^{4}+fxy^{2}+gy^{3}+hx^{3}y+jx^{2}y^{2}+kxy^{3}+ly^{4})=0
$

conic with one branch thru origin, quartic with two branches thru origin,
all three branches sharing the same tangent line.

\bigskip
\begin{figure}[!h]
  \caption{a=2; b=5/2, 3, 7/2, 4, 5, or 6.}
  \begin{center}
{\includegraphics{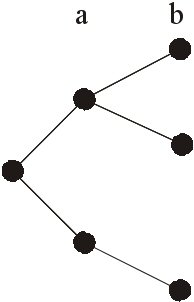}}
\end{center}
\end{figure}
\begin{figure}[!h]
  \caption{a=2, 3, 4}
  \begin{center}
{\includegraphics{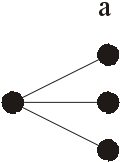}}
\end{center}
\end{figure}
6. $%
(y+ax^{2}+bxy+cy^{2}+dx^{3}+ex^{2}y+fxy^{2}+gy^{3})(y+hx^{2}+jxy+ky^{2}+lx^{3}+mx^{2}y+nxy^{2}+py^{3})
$

two cubics each with one branch thru origin and a tangent line in common.
\begin{figure}[!h]
  \caption{a=2, 3, 4, 5, 6, 7, 8, or 9.}
  \begin{center}
{\includegraphics{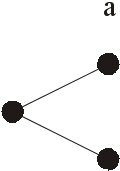}}
\end{center}
\end{figure}

7.$(y+ax^{2}+bxy+cy^{2}+dx^{3}+ex^{2}y+fxy^{2}+gy^{3})(y^{2}+hxy+jx^{3}+kx^{2}y+lxy^{2}+my^{3})$
cubic with one branch thru origin, cubic with two branches thru
origin, one tangent line in common.
\bigskip
\begin{figure}[!h]
  \caption{a=1; b=2, 3, 4, 5, 6, 7 or 8.}
  \begin{center}
{\includegraphics{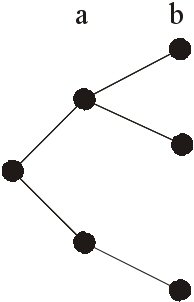}}
\end{center}
\end{figure}
\clearpage
\bigskip

\bigskip

\bigskip

8. $%
(y^{2}+axy+bx^{3}+cx^{2}y+dxy^{2}+ey^{3})(y^{2}+fxy+gx^{3}+hx^{2}y+jxy^{2}+ky^{3})=0
$

two cubics, each with two branches thru origin, with one or two tangent
lines in common.
\begin{figure}[!h]
  \caption{a=1; b=2, 3, 4, 5, or 6.}
  \begin{center}
{\includegraphics{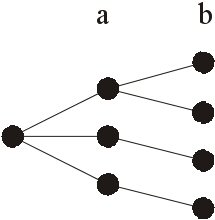}}
\end{center}
\end{figure}
\begin{figure}[!h]
  \caption{a=1; b=2; c= 3, 4, or 5 and a=1;
b=3; c=4.}
  \begin{center}
{\includegraphics{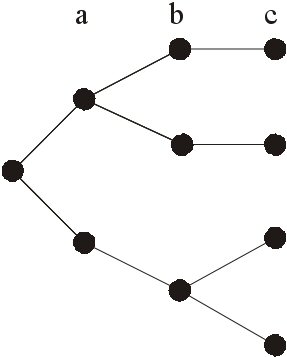}}
\end{center}
\end{figure}

\bigskip
\begin{figure}[!h]
  \caption{a=1; b=2 or 3.}
  \begin{center}
{\includegraphics{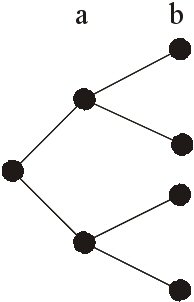}}
\end{center}
\end{figure}
9. \ $%
(y^{2}+bx^{3}+cx^{2}y+dxy^{2}+ey^{3})(y^{2}+gx^{3}+hx^{2}y+jxy^{2}+ky^{3})=0$

two cubics, each with two branches thru origin, with a double tangent line
in common.
\begin{figure}[!h]
  \caption{a=3/2.}
  \begin{center}
{\includegraphics{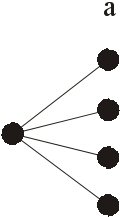}}
\end{center}
\end{figure}
\clearpage
\begin{figure}[!h]
  \caption{a=3/2; b=2 or 5/2.}
  \begin{center}
{\includegraphics{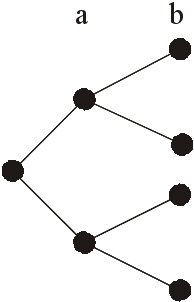}}
\end{center}
\end{figure}
\bigskip

\bigskip

\section{Summary of Singular Points of Reducible Sextic Curves}

\bigskip

\subsection{Multiplicity 2}

\bigskip

\bigskip

\bigskip

\begin{figure}[!h]
  \caption{a=1, 3/2, 2, 5/2, 3, 7/2, 4, 9/2, 5, 11/2, 6, 13/2, 7, 8, 9, 10.}
  \begin{center}
{\includegraphics{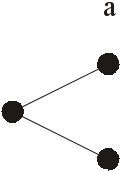}}
\end{center}
\end{figure}
\subsection{Multiplicity 3}

\bigskip

\begin{figure}[!h]
  \caption{a=1, 4/3, 3/2, 5/3, 2, 3, 4.}
  \begin{center}
{\includegraphics{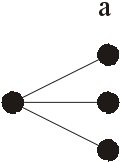}}
\end{center}
\end{figure}
\clearpage
\bigskip
\begin{figure}[!h]
  \caption{a=1; b=3/2, 2, 5/2, 3, 7/2, 4, 9/2, 5, 11/2, 6, 13/2, 7, 8 and a=2; b=5/2, 3, 7/2, 4, 9/2, 5,
11/2, 6, 13/2 and a=3; b=4.}
  \begin{center}
{\includegraphics{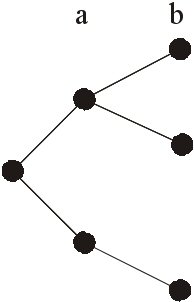}}
\end{center}
\end{figure}
\bigskip

\subsection{Multiplicity 4}

\bigskip
\begin{figure}[!h]
  \caption{a=1, 5/4, 4/3, 3/2.}
  \begin{center}
{\includegraphics{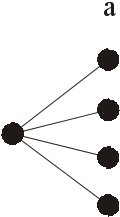}}
\end{center}
\end{figure}
\bigskip
\begin{figure}[!h]
  \caption{a=1; b=4/3, 3/2, 5/3, 2.}
  \begin{center}
{\includegraphics{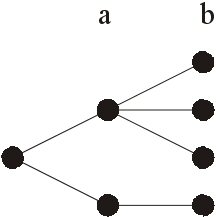}}
\end{center}
\end{figure}
\bigskip
\begin{figure}[!h]
  \caption{a=1; b=3/2, 2, 3 and a=3/2; b=2,
5/2.}
  \begin{center}
{\includegraphics{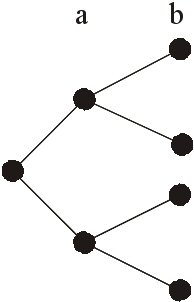}}
\end{center}
\end{figure}
\clearpage
\bigskip
\begin{figure}[!h]
  \caption{a=1; b=3/2, 2, 5/2, 3, 7/2, 4,
9/2, 5, 6 and a=3/2; b=2.}
  \begin{center}
{\includegraphics{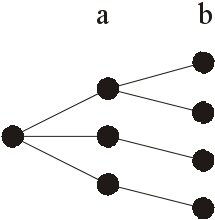}}
\end{center}
\end{figure}
\bigskip
\begin{figure}[!h]
  \caption{a=1; b=3/2; c=2,3,4,5,6 and b=2;
c=5/2, 3, 7/2, 4, 9/2, 5 and b=5/2; c=3 and b=3; c=7/2, 4, 9/2.}
  \begin{center}
{\includegraphics{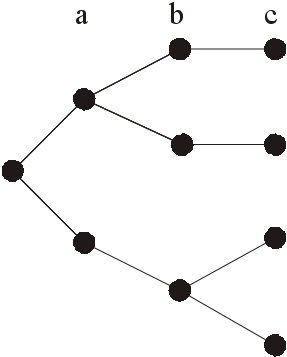}}
\end{center}
\end{figure}
\bigskip
\begin{figure}[!h]
  \caption{a=1; b=2; c=5/2, 3, 7/2, 4, 9/2, 5.}
  \begin{center}
{\includegraphics{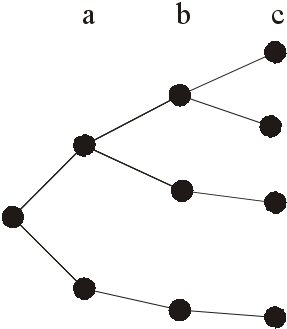}}
\end{center}
\end{figure}
\subsection{Multiplicity 5}

\bigskip

\bigskip

\bigskip

\bigskip
\begin{figure}[!h]
  \caption{a=1, 5/4.}
  \begin{center}
{\includegraphics{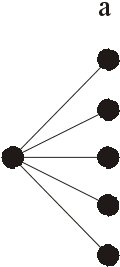}}
\end{center}
\end{figure}
\clearpage

\bigskip
\begin{figure}[!h]
  \caption{a=1; b=3/2, 2.}
  \begin{center}
{\includegraphics{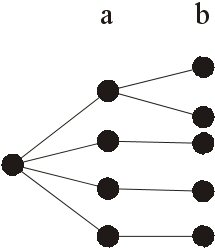}}
\end{center}
\end{figure}
\bigskip

\bigskip
\begin{figure}[!h]
  \caption{a=1; b=4/3, 5/4.}
  \begin{center}
{\includegraphics{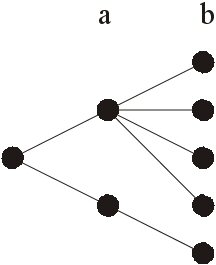}}
\end{center}
\end{figure}
\bigskip

\bigskip
\begin{figure}[!h]
  \caption{a=1; b=3/2, 2.}
  \begin{center}
{\includegraphics{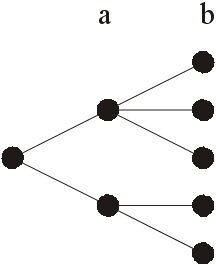}}
\end{center}
\end{figure}
\bigskip

\bigskip
\begin{figure}[!h]
  \caption{a=1; b=4/3, 3/2, 2.}
  \begin{center}
{\includegraphics{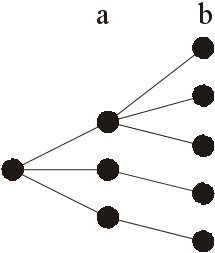}}
\end{center}
\end{figure}
\bigskip

\bigskip
\begin{figure}[!h]
  \caption{a=1; b=3/2.}
  \begin{center}
{\includegraphics{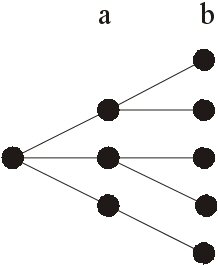}}
\end{center}
\end{figure}
\bigskip

\bigskip
\begin{figure}[!h]
  \caption{a=1; b=4/3, 3/2; c=2.}
  \begin{center}
{\includegraphics{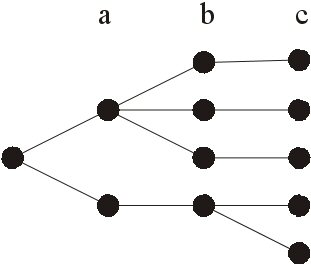}}
\end{center}
\end{figure}
\bigskip

\bigskip
\begin{figure}[!h]
  \caption{a=1; b=3/2; c=2.}
  \begin{center}
{\includegraphics{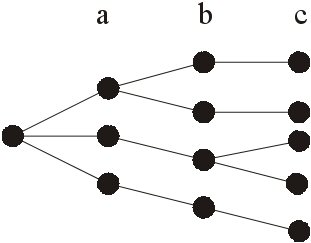}}
\end{center}
\end{figure}
\bigskip

\bigskip

\bigskip

\bigskip

\subsection{\protect\bigskip Multiplicity 6}

\bigskip
\begin{figure}[!h]
  \caption{a=1.}
  \begin{center}
{\includegraphics{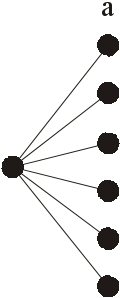}}
\end{center}
\end{figure}

\textsc{David A. Weinberg: Department of Mathematics and Statistics, Texas
Tech University, Lubbock, TX 79409-1042}

e-mail address: david.weinberg@ttu.edu

\bigskip

\textsc{Nicholas J. Willis: Department of Mathematics and Computer Science,
Whitworth University, Spokane, WA 99251}

e-mail address: nwillis@whitworth.edu\textsc{\ }

\end{document}